\documentclass[11pt,leqno]{article}

\usepackage[french]{babel}
\usepackage[ansinew]{inputenc}
\usepackage{enumerate}
\usepackage{amsfonts}
\usepackage{amssymb}
\usepackage{amsmath}
\usepackage{amscd}
\usepackage{xypic}
\newtheorem{theo}{Th\'eor\`eme}[section]
\newtheorem{prop}[theo]{Proposition}
\newtheorem{lem}[theo]{Lemme}
\newtheorem{cor}[theo]{Corollaire}
\newtheorem{rem}[theo]{Remarque}

\def \endproof{\hfill$\square$ \bigskip}
\def\id{\textrm{\scriptsize id}}
\def\o{\textrm{O}}
\def\K{\mathbb{K}}
\def\D1{
\begin{equation*}
\raisebox{1.6cm}
{\xymatrix{ 
\ar[r]^-\id         & (\K,\o) \ar[r]^-\id \ar[d]^{\Gamma_i} & (\K,\o) \ar[r]^-\id \ar[d]^{\Gamma_{i-1}} & \cdots \ar[r]^-\id     & (\K,\o) \ar[r]^-\id \ar[d]^{\Gamma_1} & (\K,\o) \ar[d]^{\Gamma_0} \\
\ar[r]^-{\Pi_{i+1}} & (E_i,\o_i) \ar[r]^-{\Pi_i}            & (E_{i-1},\o_{i-1}) \ar[r]^-{\Pi_{i-1}}    & \cdots \ar[r]^-{\Pi_2} & (E_1,\o_1) \ar[r]^-{\Pi_1}            & (E_0,\o_0)                \\
\ar[r]^-{\omega_{i+1}}   & (Z_i,\o_i) \ar@{^{(}->}[u]_{J_i} \ar[r]^-{\omega_i} & (Z_{i-1},\o_{i-1}) \ar@{^{(}->}[u]_{J_{i-1}} \ar[r]^-{\omega_{i-1}}      & \cdots \ar[r]^-{\omega_2}   & (Z_1,\o_1) \ar@{^{(}->}[u]_{J_1} \ar[r]^-{\omega_1}              & (Z_0,\o_0) \ar@{^{(}->}[u]_{J_0}
}}
\end{equation*}
}

\begin{document}
\baselineskip=16pt
\title{Fonction asymptotique de Samuel des sections hyperplanes et multiplicité}
\author{M. Hickel}
\date{ }\maketitle
\renewcommand{\abstractname}{Abstract}
\begin{abstract}
\baselineskip=16pt
 Let $(A,\mathfrak{m}_A,k)$ be a local noetherian ring and $I$ an $\mathfrak{m}_A$-primary ideal. The asymptotic Samuel function (with respect to  $I$) $\overline{v_I}$: $A\longrightarrow \mathbb{R}\cup \{+\infty\}$ is defined by $\overline{v_I}(x)=lim_{k\rightarrow +\infty}\frac{ord_I(x^k)}{k}$, $\forall x \in A$. Similary, one defines for another ideal $J$, $\overline{v_I}(J)$ as the minimum of $\overline{v_I}(x)$ as $x$ varies in $J$. Of special interest is the rational number $\overline{v_I}(\mathfrak{m}_A)$. We study the behavior of the asymptotic Samuel function (with respect to $I$) when passing to hyperplanes sections of $A$ as one does for the theory of mixed multiplicities. 
\footnote{2000 \it Mathematics Subject Classification; Primary 13B22; Secondary 13C15, 13F25.
 Key words and phrases: Asymptotic Samuel function, Hyperplanes sections, Integral Closure of ideals, Multiplicity, \L ojasiewicz inequalities.}\end{abstract}

 \section{Introduction}
 Soient $(A,\mathfrak{m}_A,k)$ un anneau local noetherien et $I$ un id\'eal $\mathfrak{m}_A$-primaire. La fonction asymptotique de Samuel (par rapport \`a $I$) est la fonction $\overline{v_I}$ : $A\longrightarrow \mathbb{R}_{\geq 0}\cup \{+\infty\}$ d\'efinie par :
$$\overline{v_I}(x)=Lim_{k \rightarrow + \infty}\frac{ord_I(x^k)}{k}$$
o\`u $ord_I(x)=Max\{m\in \mathbb{N}/ x \in I^m\}$ si $x\neq 0$ et $ord_I(0)=+\infty$. De fa\c{c}on similaire, on d\'efinit pour un id\'eal $J$ de $A$ le nombre $\overline{v_I}(J)$ par : 
$$ \overline{v_I}(J)=Min_{x \in J}\overline{v_I}(x).$$
La fonction asymptotique de Samuel prend en fait des valeurs rationnelles et peut se calculer \`a l'aide des \textit {valuations de Rees} de $I$. Elle joue un r\^ole primordial dans l'\'etude des cl\^otures int\'egrales des puissances de $I$. Nous renvoyons \`a [H-S], [L-T], [R1] pour les r\'esultats et les compl\'ements essentiels concernant ces notions. En particulier le nombre rationnel $\overline{v_I}(\mathfrak{m})$, ou son inverse $\nu _I(\mathfrak{m})$=$\overline{v_I}(\mathfrak{m})^{-1}$ que nous appellerons \textit{exposant de \L ojasiewicz de $I$} (voir § 2) joue un r\^ole important dans l'\'etude de la topologie des singularit\'es c.f. [T1,2]. De la m\^eme mani\`ere que la notion de \textit{multiplicit\'e mixte} permet en particulier de rendre compte du comportement de la $I$-multiplicit\'e de $A$ apr\`es avoir intersect\'e celui-ci par des sections hyperplanes suffisamment g\'en\'eriques, nous avons cherch\'e \`a \'etablir des r\'esultats similaires pour la $I$-fonction asymptotique de Samuel apr\`es sections hyperplanes g\'en\'eriques. Le r\'esultat principal du pr\'esent travail est le suivant.
\begin{theo}\textrm{ }\\
Soit $(A,\mathfrak{m},k)$ un anneau local r\'egulier d'\'egale caract\'eristique z\'ero, dont on notera $n+1$ la dimension. Consid\'erons un id\'eal $I\subset A$, $\mathfrak{m}$-primaire, de multiplicit\'e $e(I)$ et d'exposant de \L ojasiewicz $\nu _I(\mathfrak{m})=\overline{v_I}(\mathfrak{m})^{-1}=\nu _I^{(n+1)}$.
\begin{itemize}
\item[1)] Pour tout $i$, $1\leq i\leq n$, il existe un ouvert de Zariski dense $U^{(i)}=U^{(i)}(I)$  $\subset  G_k(i,n+1)$ de la Grassmanienne des $i$-plans de $k^{n+1}$ tel que $\forall H \in U^{(i)}$, $H$ d\'efini par l'annulation des $n+1-i$ formes lin\'eaires 
$$h_j(X_0,X_1,\ldots,X_n)=\sum _{0\leq k \leq n}a_{k,j}.X_j,$$
le nombre rationnel $\overline{v_{I.A_H}}(\mathfrak{m}_H)$, o\`u $A_H=A/(h_1,\ldots, h_{n+1-i})$ et $\mathfrak{m}_H$ est l'id\'eal maximal de $A_H$, est ind\'ependant de $H\in U^{(i)}$. Le nombre rationnel $\nu ^{(i)}_I=(\overline{v_{I.A_H}}(\mathfrak{m}_H))^{-1}$ (ind\'ependant de $H \in U^{i}$) est appel\'e le i\`eme exposant de \L ojasiewicz de $I$ ou encore l'exposant de \L ojasiewicz de I \textrm{restreint \`a un i-plan g\'en\'erique}.
\item[2)] On a  :
$$e(I)\leq \nu _I^{(1)}\times\nu ^{(2)}_I\times \ldots.\times \nu _I^{(n)}\times \nu _I^{(n+1)}.$$
\end{itemize}
\end{theo} 
Notons que $\nu ^{(1)}_I\leq\nu ^{(2)}_I\leq\ldots \leq \nu ^{(n)}_I\leq\nu ^{(n+1)}_I$ et que $\nu ^{(1)}_I$ n'est autre que $ord_{\mathfrak{m}}(I)$. Un cas important parmi les cas d'\'egalit\'es dans 2) est fourni par les id\'eaux $I$ tels qu'il existe des entiers $b\in \mathbb{N}^*$, $1\leq a_1\leq a_2\ldots \leq a_{n+1}$ tels que 
$$\overline{I^b}=\overline{(l_1^{a_1},l_2^{a_2},\ldots,l_{n+1}^{a_{n+1}})},$$
o\`u $l_1,\ldots,l_{n+1}$ est un syst\`eme r\'egulier de param\`etres de $A$. En effet, dans ce cas $\nu _I^{(i)}=\frac{a_i}{b}$ et $b^{n+1}e(I)=\Pi _{1\leq i \leq n+1}a_i$. Une discussion et une description plus compl\`ete des cas d'\'egalit\'es de $2)$ est donn\'ee plus bas (c.f  5)).\\
Quelques mots sur la preuve du r\'esultat ci-dessus. Tout d'abord, il ne nous a pas paru ais\'e de d\'ecrire la variation des valuations de Rees de $I.A_H$ lorsque $H$ parcourt la Grassmanienne et donc de suivre par ce biais la variation de $\overline{v_{I.A_{H}}}(\mathfrak{m}_H)$. Ainsi la preuve du premier point de 1.1 proc\`ede de mani\`ere indirecte. On commence par montrer un r\'esultat ind\'ependant (Th 3.1) qui nous semble pouvoir pr\'esenter un inter\^et par lui-m\^eme. Si $A$ est local r\'egulier d'\'egale caract\'eristique z\'ero, nous montrons comment on peut calculer $\overline{v_I}(\mathfrak{m})$ \`a l'aide de certains polyn\^omes caract\'eristiques canoniquement associ\'es \`a $I$. Ce proc\'ed\'e, ind\'ependant de la connaissance des valuations de Rees de $I$, fait l'objet de la section 3).  Nous rappellons au pr\'ealable (section 2) quelques r\'esultats et notions que nous utiliserons ensuite. C'est dans la preuve de 3.1 qu'interviennent les hypoth\`eses de r\'egularit\'e et d'\'egale caract\'eristique z\'ero. On utilise ensuite l'existence de r\'eductions jointes  et des techniques de bases standards (dans le sens de Grauert-Hironaka [A-H-V], [H]) et de variations de diagrammes des exposants initiaux telles que d\'evelopp\'ees par Bierstone-Milman [B-M] pour montrer que les polyn\^omes caract\'eristiques en question varient de fa\c{c}on agr\'eable lorsque $H$ d\'ecrit un ouvert de Zariski convenable de la Grassmanienne. Ceci est expos\'e \`a la section 4) et permet d'obtenir le premier point de 1.1. La majoration de la multiplicit\'e $e(I)$ est obtenue ensuite comme cons\'equence de cela et d'une g\'en\'eralisation de la loi d'associativit\'e pour les multiplicit\'es que l'on peut trouver dans le livre de D.G. Northcott [N] (Nous ignorons si ce r\'esultat lui est d\^u). Nous d\'ecrivons ensuite certains cas d'\'egalit\'e dans l'in\'egalit\'e 2) \`a la section 5).
\section{Rappels et notations}
\subsection{L'algorithme de division formelle de Grauert-Hironaka}
Soit $\alpha=(\alpha_{1},\ldots,\alpha_{n})\in\mathbb{N}^{n}$, on  notera: $\vert\alpha\vert=\alpha_{1}+\cdots+\alpha_{n}$. $\mathbb{N}^{n}$ est totalement ordonn\'e par l'ordre lexicographique sur les $n+1$ upplets  $(\vert\alpha\vert,\alpha_{1},\ldots,\alpha_{n})$. Si $A$ est un anneau commutatif unitaire intègre, et $f$ un élément non nul de $A\lbrack\lbrack X \rbrack\rbrack=A\lbrack\lbrack X_{1},\ldots,X_{n}\rbrack\rbrack$, nous noterons par $\nu (f)$ son exposant initial. C'est \`a dire si : 
$$ f=\sum_{\alpha\in\mathbb{N}^n}a_{\alpha}.X^{\alpha}, \hspace{3pt} Supp(f)=\{\alpha\in\mathbb{N}^n \vert a_{\alpha}\not=0\},\hspace{3pt} \nu (f)=Min\{\alpha\in\mathbb{N}^n \vert a_{\alpha}\not=0\}.$$ 
De m\^eme, on notera $Init(f)$ le  coefficient du mon\^ome initial de $f$ i.e. $Init(f)=a_{\nu (f)}$.
Pour un id\'eal, $I \subset A[[X]]$, on notera $\bigtriangleup _{I}$ le diagramme des exposants initiaux de $I$, i.e.: 
$$\bigtriangleup _{I}=\{\alpha \in \mathbb{N}^{n} \vert \exists g \in I \textit{ tel que }  \nu(g)=\alpha\}.$$
On a : $\bigtriangleup _{I} + \mathbb{N}^{n} =\bigtriangleup _I$. Si $N \subset \mathbb{N}^{n}$ satisfait $N + \mathbb{N}^{n} =N$ (i.e. est stable par translation), le lemme de Dickson nous assure de l'existence d'une unique partie finie de $\mathbb{N}^{n}$, $\{\alpha^{1},\ldots,\alpha^{p}\}$, telle que : 
$$N=\cup_{1\leq i\leq p}(\alpha^{i} +  \mathbb{N}^{n}) \textit { et } \alpha ^{j} \notin \cup_{i\neq j} (\alpha^{i} +  \mathbb{N}^{n})$$
Les $\alpha^{i}$ sont dits les sommets de $N$. L'ensemble $\mathcal{D}(n)$ des parties  de $\mathbb{N}^{n}$ stables par translation (i.e. satisfaisant $N +\mathbb{N}^{n}=N$) est totalement ordonn\'e comme suit. Soient $N_{1},N_{2}\in \mathcal{D}(n)$. Pour chaque $i=1,2$, soient $\beta ^{k}_{i}$, $k=1,\ldots,t_{i}$ les sommets de $N_{i}$ index\'es dans l'ordre croissant. Apr\`es avoir \'eventuellement permut\'e $N_{1}$ et $N_{2}$, il existe $t\in \mathbb{N}$ tel que : $\beta^{k}_{1}=\beta^{k}_{2}$, $1\leq k\leq t$ et (1) $t_{1}=t=t_{2}$, ou bien (2) $t_{1}>t=t_{2}$, ou bien (3) $t_{1},t_{2}>t$ et $\beta_{t+1}^{1}<\beta_{t+1}^{2}$. Dans le cas (1), $N_{1}=N_{2}$. Dans les cas (2) et (3), $N_{1}<N_{2}$. Il revient au m\^eme de dire que la suite $(\beta ^1_1,\ldots,\beta ^{t_1}_1,\infty,\ldots)$ est strictement plus petite que la suite $(\beta ^1_2,\ldots,\beta ^{t_2}_2,\infty,\ldots)$ pour l'ordre lexicographique, avec la convention $\beta <\infty $ pour tout $\beta \in \mathbb{N}^n$.\\
Si $\beta ^1,\ldots,\beta ^t \in \mathbb{N}^n$, on leur associe la partition suivante de $\mathbb{N}^n$ :
$$\bigtriangleup =\cup _{1\leq i \leq t}(\beta ^i+\mathbb{N}^n),\quad \overline{\bigtriangleup }=\mathbb{N}^n-\bigtriangleup ,\quad \bigtriangleup _i=(\beta ^i+\mathbb{N}^n)-\cup _{k<i}(\beta ^k+\mathbb{N}^n)$$
Ainsi $\mathbb{N}^n$ est l'union disjointe : $\bigtriangleup \cup \overline{\bigtriangleup }=(\cup _{1\leq i \leq t}\bigtriangleup _i)\cup \overline{\bigtriangleup }$. Nous utiliserons le r\'esultat suivant.
\begin{theo} \textit{(Le th\'eor\`eme de division formelle de Grauert-Hironaka [A-H-V], [B-M], [G])}\\
Soient $A$ un anneau (commutatutif unitaire) int\`egre et $F_1,\ldots,F_t$ des \'el\'ements non nuls de $A[[X_1,\ldots,X_n]]$. Notons  $\beta ^i=\nu (F_i)$ l'exposant initial de $F_i$, $a_i=Init(F_i)$ le coefficient du mon\^ome initial de $F_i$. Soit $S$ la partie multiplicativement ferm\'ee de $A$ engendr\'ee par les $a_i$ i.e. $S=\{a_1^{m_1},\ldots,a_t^{m_t}, \quad m_i\in \mathbb{N}\}$. Alors tout \'el\'ement $F$ de $A[[X]]$ (ou de $S^{-1}.A[[X]]$) peut s'\'ecrire de mani\`ere unique sous la forme :
$$F=\sum_{1\leq i \leq t}D_iF_i+R,\quad D_i,R \in S^{-1}.A[[X]]$$
avec $Supp(R)\subset \overline{\bigtriangleup }$ et $Supp(D_i)+\beta _i\subset \bigtriangleup _i$, o\`u les $\bigtriangleup _i$ et $\overline{\bigtriangleup }$ sont d\'efinis comme ci-dessus.
\end{theo}
\subsection{Brefs rappels sur la fonction asymptotique de Samuel et les valuations de Rees}
Soient  $(A,\mathfrak{m}_A,k)$ un anneau local noetherien int\`egre et $I$ un id\'eal de $A$. Comme nous l'avons rappel\'e plus haut la fonction asymptotique de Samuel relativement \`a $I$ est d\'efinie par :
$$ \forall x \in A,\quad \overline{v_I}(x)=Lim_{k \rightarrow + \infty}\frac{ord_I(x^k)}{k}.$$
Nous rappelons ici tr\`es bri\`evement comment se calcule $\overline{v_I}$ \`a l'aide de valuations discr\`etes de rang $1$ positives sur $A$. Pour de plus amples compl\'ements nous renvoyons \`a [H-S] chap.6 et 10. Par la terminologie \guillemotleft $v$ une valuation discr\`ete de rang $1$ positive sur $A$ \guillemotright nous entendrons la donn\'ee d'une valuation $v$ associ\'ee \`a un anneau de valuation discr\`ete de rang $1$, $A_v$, entre $A$ et son corps des fractions $K=Frac(A)$ (i.e. $A\subset A_v \subset K$) tel que $\mathfrak{m}_v \cap A=\mathfrak{m}_A$. D\'esignons par $\Lambda _A$ l'ensemble des telles valuations. Alors :
$$ \forall x \in A,\quad \overline{v_I}(x)=Inf_{v \in \Lambda _A}\frac{v(x)}{v(I)}=Min_{v \in \Lambda _A}\frac{v(x)}{v(I)}$$
De plus, il existe un ensemble $\mathcal{R_I}$ fini non redondant de telles valuations, unique \`a l'\'equivalence pr\`es des valuations, telles que :
$$ \forall x \in A,\quad \overline{v_I}(x)=Min_{v \in \mathcal{R}_I}\frac{v(x)}{v(I)}$$
Ces valuations particuli\`eres sont appel\'ees l'ensemble des valuations de Rees de $I$. Diverses constructions des valuations de Rees de $I$ sont expos\'ees dans [H-S] chap. 10. Un point de vue plus g\'eom\'etrique dans le cadre de la g\'eom\'etrie analytique complexe est dans [L-T]. En particulier [L-T], si $A=\mathcal{O}_{X,x}$ est la fibre en $x$ du faisceau structural d'un espace analytique complexe $X$ et $I=(f_1,\ldots,f_d)$ est un id\'eal $\mathfrak{m}_{X,x}$-primaire, le nombre $\overline{v_I}(\mathfrak{m}_{X,x})^{-1}$ est le plus petit nombre r\'eel positif $\alpha $ tel qu' il existe un voisinage de $V_x$ de $x$ dans $X$ et  une constante $C$ tels que :
$$\forall z \in V_x,\quad \sum_{i=1}^d \vert f_i(z)\vert \geq C\vert\vert z-x\vert \vert^{\alpha},
 $$
o\`u $\vert \vert \quad \vert \vert$ d\'esigne une norme arbitraire sur $\mathbb{C}^n$ dans lequel on a plong\'e $(X,x)$.
Ceci est un cas particulier \textit{d'in\'egalit\'es de \L ojasiewicz}. Pour cette raison et dans un contexte g\'en\'eral nous appellerons le nombre $\overline{v_I}(\mathfrak{m}_A)^{-1}$ l'exposant de \L ojasiewicz de $I$.

\section{Polyn\^omes caract\'eristiques et calcul de la fonction asymptotique de Samuel}

Soit $(A,\mathfrak{m})$ un anneau local r\'egulier d'\'egale caract\'eristique z\'ero, de dimension $n$. Soit $I$ un id\'eal $\mathfrak{m}$-primaire et $g$ un \'el\'ement de $A$. Nous cherchons un proc\'ed\'e qui permet de calculer $\overline{v}_I(g)$ sans recours \`a la description des valuations de Rees, en particulier sans recours \`a l'\'eclatement normalis\'e de centre $I$ (i.e. sans recours \`a la cl\^oture int\'egrale de l'alg\`ebre de Rees de $I$). Nous allons voir que pour tout $g \in A$, il y a une  relation de d\'ependance int\'egrale de degr\'e $e(I)$, canonique, qui calcule $\overline{v}_I(g)$. Soit $\hat{A}$ le compl\'et\'e pour la topologie $\mathfrak{m}$-adique de $A$. Par fid\`ele platitude $I^p.\hat{A}\cap A=I^p$. Il en d\'ecoule alors que pour tout $g\in A$, $\overline{v}_I(g)=\overline{v}_{I.\hat{A}}(g)$. On peut donc supposer que $A$ est complet et donc par le th\'eor\`eme de structure de Cohen supposer que $A$ est isomorphe \`a $k[[X_1,\ldots,X_n]]$ o\`u $k\simeq A/\mathfrak{m}$. Supposons dans un premier temps que $I=(f_1,\ldots,f_n)$ est engendr\'e par une suite r\'eguli\`ere. La construction de base est alors la suivante. Notons $F^{*}$ le morphisme de $A_1=k[[U_1,\ldots,U_n]]$ dans $A_2=k[[X_1,\ldots,X_n]]$ d\'efini par $U_i\longrightarrow f_i(X)$. Le morphisme $F^{*}$ est quasi-fini car $dim_k(A_2/\mathfrak{m}_1.A_2)<+\infty$, il est donc fini par le th\'eor\`eme de pr\'eparation formel c.f. [La] chap. 8. Il en r\'esulte aussit\^ot qu'il est injectif, puisque $A_2$ est entier sur $A_1/Ker(F^*)$, ces deux anneaux ont donc la m\^eme dimension et par suite $Ker(F^*)=0$. Maintenant le crit\`ere local de platitude (c.f. [M]), nous assure que $A_2$ est un $A_1$ module plat de type fini car $(f_1,\ldots,f_n)$ est une suite r\'eguli\`ere de $A_2$. Par cons\'equent les anneaux \'etant locaux, $A_2$ est un $A_1$-module libre de type fini, et son rang n'est autre que $dim_k (A_2/\mathfrak{m_1}.A_2)=dim_k (A_2/I)=e(I)$, puisque $A_2$ est de Cohen-Macauley.  Si $g$ est un \'el\'ement de $A_2$, l'op\'erateur de multiplication par g,  $m_g$: $A_2\longrightarrow A_2$ a donc un polyn\^ome caract\'eristique comme \'el\'ement de $End_{A_1}(A_2)$. Notons $P_g(Y,U) \in A_1[Y]$ ce polyn\^ome caract\'eristique :
$$P_g(Y,U)=Y^{e(I)}+\sum_{k=1}^{e(I)}a_i(U_1,\ldots,U_n)Y^{e(I)-i} \in A_1[Y].$$ Clairement, par le th\'eor\`eme d'Hamilton-Cayley on a : $P_g(g,f_1,\ldots,f_n)=0$.  Cette relation de d\'ependance int\'egrale particuli\`ere calcule toujours $\overline{v}_I(g)$. On a en effet le r\'esultat suivant : 
\begin{theo}\textrm{ }\\
Soient $I=(f_1,\ldots,f_n)\subset A_2$ et $g \in A_2$ comme ci-dessus, et $P_g(Y,U)\in A_1[Y]$ son polyn\^ome caract\'eristique. Posons :
$$ \alpha =\frac{p}{q}=Min_{1\leq i \leq e(I)}( \frac{Ord_U a_i(U)}{i}),\textit{ o\`u } ord_U a=Max\{k \in \mathbb{N}/a \in (U)^k\}$$
Alors : $\overline{v}_I(g)=\alpha =\frac{p}{q}$.
\end{theo} 
\noindent \textit{Preuve :}\\
Nous constatons d'abord que $\overline{v}_I(g)\geq \frac{p}{q}$. Notons $K_2=Frac(A_2)$. En effet, soit $v$ une valuation discr\`ete $K_2\longrightarrow \mathbb{Z}$ positive sur $A_2$ de rang $1$. Puisque :
$$g^{e(I)}+\sum_{i=1}^{e(I)}a_i(f_1,\ldots,f_n)g^{e(I)-i}=0.$$
On a :  $$(*)\hspace{10pt}e(I)v(g)\geq Min_{1\leq i \leq e(I)}(v(a_{i}(f_1,\ldots,f_n))+(e(I)-i)v(g)).$$ 
Soit $i_0$ r\'ealisant le minimum dans  le membre de droite de cette in\'egalit\'e. Par d\'efinition on a : $ord_{(U)} a_{i_0}(U)\geq i_0\frac{p}{q}$. Donc $ord_I a_{i_0}(f_1,\ldots,f_n)\geq i_o\frac{p}{q}$. Par cons\'equent : $v(a_{i_0}(f_1,\ldots,f_n))\geq i_0\frac{p}{q}v(I)$. Il vient ainsi simplifiant l'in\'egalit\'e $(*)$ par $(e(I)-i_0)v(g)$ : $i_0 v(g)\geq i_0. \frac{p}{q}.v(I)$, o\`u encore $v(g)\geq \frac{p}{q}. v(I)$. D'o\`u l'on d\'eduit que $\overline{v}_I(g) \geq \frac{p}{q}$, puisque $\overline{v}_I(g)$ se calcule comme le minimum des $v(g)/v(I)$ lorsque $v$ parcourt l'ensemble des valuations de Rees de $I$. Supposons maintenant que l'in\'egalit\'e $\overline{v}_I(g)\geq \frac{p}{q}$ soit stricte i.e.  $\overline{v}_I(g) >\frac{p}{q}$. Nous allons voir que l'on aboutit \`a une contradiction.  Remarquons que pour cela, on peut supposer que $k$ est alg\'ebriquement clos. En effet, sinon soit $\overline{k}$ une cl\^oture alg\'ebrique de $k$. Le morphisme $k[[X_1,\ldots,X_n]]\longrightarrow \overline{k}[[X_1,\ldots,X_n]]$ \'etant fid\`element plat, on a $\overline{v}_I(g)=\overline{v}_{I.\overline{k}[[X]]}(g)$, et de m\^eme les polyn\^omes caract\'eristiques de $g$  consid\'er\'e comme \'el\'ement de $k[[X]]$ ou comme \'el\'ement de $\overline{k}[[X]]$ coincident. On supposera donc que $k$ est alg\'ebriquement clos. Soit $K_1=Frac(A_1)$ et notons toujours $F^{*}$: $K_1\longrightarrow K_2$ le morphisme injectif induit par $F^{*}$: $A_1\longrightarrow A_2$. Comme $A_2$ est un $A_1$-module libre de type fini, on a facilement que pour tout $b\in A_2-(0)$, $1/b \in K_1.A_2$ et donc $K_2$ est une extension alg\'ebrique finie de $K_1$ et $[K_2:K_1]=e(I)$. Comme l'on a suppos\'e que $A$ \'etait d'\'egale caract\'eristique z\'ero, cette extension est s\'eparable et quitte \`a faire un changement de variables lin\'eaires sur les $X_i$, par le th\'eor\`eme de l'\'el\'ement primitif on peut supposer que cette extension est engendr\'ee par $X_1$ i.e. $K_2=K_1(X_1)$.  Maintenant soit $\frac{p'}{q'}$ tel que $\overline{v}_I(g)>\frac{p'}{q'}>\frac{p}{q}$. Puisque $\overline{v}_I(g)=Lim_{m\rightarrow \infty}\frac{ord_I(g^m)}{m}$, il existerait $m_0$ tel que pour $m\geq m_0$ on ait $ord_I(g^m)\geq m \frac{p'}{q'}$. D\'esignant par $[\quad ]$ la partie enti\`ere sup\'erieure, on aurait donc $g^m \in I^{[m\frac{p'}{q'}]}$. Par suite pour tout arc $\varphi ^*$: $k[[X_1,\ldots,X_n]]\longrightarrow k[[t]]$, on aurait : $$m.ord_t(\varphi ^*(g))\geq [m.\frac{p'}{q'}]Min_{1\leq i\leq n}ord_t(\varphi ^*(f_i))$$ 
et donc :
$$(\bullet)\hspace{5pt}ord_t(\varphi ^*(g))\geq \frac{p'}{q'}Min_{1\leq i\leq n}ord_t (\varphi ^*(f_i))>\frac{p}{q}Min_{1\leq i\leq n}ord_t (\varphi ^*(f_i)).$$ Nous allons construire un arc $\varphi ^*$ qui ne satisfait pas cette in\'egalit\'e et obtenir ainsi une contradiction, ce qui prouvera $\overline{v}_I(g)=\frac{p}{q}$. Pour cela, notons $P$ le polyn\^ome minimal de $X_1$ sur $K_1$. $P$ est de degr\'e $e(I)$ et \`a coefficients dans $A_1$, i.e. $P\in A_1[Y]$. On notera $\Delta \in k[[U_1,\ldots,U_n]]$ son discriminant (qui est non nul puisque l'extension est s\'eparable). Notons $D$ le produit :
$$D=\Delta .\Pi _{1\leq i\leq e(I)} a_i \in k[[U_1,\ldots,U_n]]$$
o\`u les $a_i$ sont les coefficients du polyn\^ome caract\'eristique de $g$ comme d\'efini ci-dessus. Consid\'erons alors $In(D)\in k[U_1,\ldots,U_n]$ la forme initiale de $D$, c'est \`a dire le polyn\^ome homog\`ene de plus bas degr\'e dans le d\'eveloppement en somme de polyn\^omes homog\`enes de $D$. $In(D)$ est le produit des formes initiales des $a_i$ et de celle de $\Delta $. Puis choisissons un point $(b_1,\ldots,b_n)\in k^n$ tel que $In(D)(b_1,\ldots,b_n)\neq 0$ et $b_i\neq 0$, $1\leq i \leq n$. On consid\`ere alors l'id\'eal $\mathfrak{a}$ de $A_2$ engendr\'e par $(b_1f_2-b_2f_1,b_1f_3-b_3f_1\ldots,b_1f_n-b_nf_1)$. Puisque le morphisme $F^{*}$ : $A_1\longrightarrow A_2$ est entier on a $haut(\mathfrak{a})=n-1$ par les th\'eor\`emes de Cohen-Seidenberg. Soit $\mathfrak{p}$ un id\'eal premier minimal de hauteur $n-1$ parmi ceux contenant $\mathfrak{a}$, soit $B$ la cl\^oture int\'egrale de $k[[X]]/\mathfrak{p}$ et $\hat{B}$ son compl\'et\'e. Alors par le th\'eor\`eme de structure de Cohen $\hat{B}$ est isomorphe \`a $k[[t]]$. On a donc un morphisme non nul (un arc) $\varphi ^{^1*}$ :
$$\varphi ^{1*}: k[[X]]\longrightarrow k[[X]]/\mathfrak{p}\longrightarrow B\longrightarrow k[[t]]$$
dont le noyau est $\mathfrak{p}$. Posons $\varphi ^1_i=\varphi ^{1*}(X_i)\in k[[t]]$. Par construction : $b_1\varphi ^{1*}(f_l)-b_l\varphi ^{1*}(f_1)=0$ et $\varphi ^{1*}(f_1)\neq 0$ (car $b_i\neq 0$ pour tout $i$ et $ker(\varphi ^{1*})=\mathfrak{p}$). Posons $u(t)=\varphi ^{1*}(f_1)/b_1$. On a donc :
$$f_l(\varphi ^1_1(t),\ldots,\varphi ^1_n(t))=b_l.u(t),\hspace{3pt}1\leq l \leq n.$$
Notant $\alpha _j$, $1\leq j\leq e(I)$ les coefficients du polyn\^ome minimal  $P$ de $X_1$ sur $K_1$. On a :
$$\varphi _1^{1}(t)^{e(I)}+\sum_{1\leq j\leq e(I)}\alpha _j(b_1.u(t),\ldots,b_n.u(t)).\varphi _1^1(t)^{e(I)-i}=0$$
Regardons l'\'equation :
$$P(t,X_1)=X_1^{e(I)}+\sum_{1\leq i \leq e(I)}\alpha _i(b_1.u(t),\ldots,b_n.u(t)).X_1^{e(I)-i}=0$$
comme une \'equation \`a coefficients dans le corps des s\'eries de puiseux $\cup _{m\geq 1}k((t^{1/m}))$ qui est alg\'ebriquement clos. $\varphi _1^1$ en est une racine,  le discriminant de ce polyn\^ome vaut $\Delta (b_1.u(t),\ldots,b_n.u(t))$ et comme $(b_1,\ldots,b_n)$ n'est pas un z\'ero de la forme initiale de $\bigtriangleup $ on a :
$$ord_t\Delta (b_1.u(t),\ldots,b_n.u(t))=ord_U(\Delta(U) ).ord_t(u(t))$$
Ainsi notre polyn\^ome $P(t,X_1)$ a $e(I)$ racines distinctes dans $\cup _{m \geq 1}k((t^{1/m}))$. Celles-ci sont toutes dans $k((t^{1/m}))$ pour $m$ convenablement choisi (assez grand). Notons $\varphi _1^1,\ldots,\varphi _1^{e(I)}$ ces racines distinctes. Elles sont en fait dans $k[[t^{1/m}]]$ puisque \'el\'ements de $k((t^{1/m}))$ et enti\`eres sur $k[[t]]$ donc sur $k[[t^{1/m}]]$. Un changement d'uniformisante $t\longrightarrow t^{1/m}$ nous donne donc $e(I)$ \'el\'ements distincts dans $k[[t]]$, que nous noterons encore $\varphi _1^1,\ldots,\varphi _1^{e(I)}$, tels qu'en posant $v(t)=u(t^m)$ on ait :
$$\varphi _1^j(t)^{e(I)}+\sum _{1\leq i \leq m}\alpha _i(b_1.v(t),\ldots,b_n.v(t)).\varphi _1^j(t)^{e(I)-i}=0,\hspace{3pt}1\leq j \leq e(I)$$
Maintenant puisque $K_2=K_1(X_1)$, on peut \'ecrire dans $K_2$  pour $s\geq 2$:
$$X_s=\sum_{m=0}^{e(I)-1}\beta _m ^s. X_1^m, \hspace{3pt}\beta _m^s\in K_1$$
Classiquement (c.f. par exemple [To] Th 7.5 p. 25), on voit que les $\beta  _m^s$ sont dans $(A_1)_{(\Delta) }$ o\`u $(\Delta)$ d\'esigne la partie multiplicativement ferm\'ee $\{1,\Delta ,\Delta ^2,\ldots\}$. On pose maintenant :
$$\varphi _s^j(t)=\sum_{0\leq m \leq e(I)-1}\beta  _m^s(b_1.v(t),\ldots,b_n.v(t)).\varphi _1^j(t)^{m} \in k((t))$$
Puisque $X_s$ est entier sur $A_1$ (et non pas simplement sur $K_1$) alors $\varphi _s^j(t)$ est \'el\'ement de $k((t))$ et entier sur $k[[t]]$, il est donc \'el\'ement de $k[[t]]$. Maintenant, on pose pour tout $j$, $1\leq j \leq e(I)$ :
$$\varphi ^j(t)=(\varphi _1^j(t),\ldots,\varphi _n^j(t))\in k[[t]]^n.$$
Par construction m\^eme on a :
$$f_l(\varphi ^j(t))=b_l.v(t), \hspace{3pt}1\leq l\leq e(I)$$
et $\varphi ^j\neq \varphi ^{j'}$ si $j\neq j'$. Nous allons voir que chacun des arcs $\varphi ^{*j}$ : $k[[X]]\longrightarrow k[[t]]$ et en particulier $\varphi^{1*}$ nous fournit une contradiction avec $(\bullet)$. 
En effet, puisque $(b_1,\ldots,b_n)$ n'annule pas les z\'eros de la forme initiale des $a_i$, on a :
\begin{align}
ord_t a_i(f_1(\varphi ^j(t)),\ldots,f_n(\varphi ^j(t))) & =ord_t a_i(b_1.v(t),\ldots,b_nv(t))\\
=ord_U a_i(U_1,\ldots,U_n)ord_tv(t)&=ord_U a_i(U_1,\ldots,U_n)Min_{1\leq l \leq n}ord \varphi ^{j*}(f_l)
\end{align}
Puisque $g(\varphi ^1(t))$ est racine de :
$$Y^{e(I)}+\sum_{1\leq i \leq m}a_i(b_1.v(t),\ldots,b_nv(t))Y^{e(I)-i}=0$$
et que les autres racines sont les $g(\varphi ^j(t))$. On a aux signes pr\`es :
$$a_i(f_1(\varphi _1^1(t),\ldots,f_n(\varphi_n ^1(t))=\sum _{1 \leq j_1 < j_2\ldots j_i\leq e(I)}g(\varphi ^{j_1}(t))\ldots g(\varphi ^{j_i}(t)).$$
Maintenant par $(\bullet)$ on a :
$$ord_t(\sum _{1 \leq j_1 < j_2\ldots j_i\leq e(I)}g(\varphi ^{j_1}(t))\ldots g(\varphi ^{j_i}(t))>i.\frac{p}{q}Min_{1\leq l \leq n}ord_t\varphi ^{1*}(f_l)$$
Donc  comparant avec $(1), (2)$ on obtient : 
$$\forall i,\quad ord_U(a_i(U_1,\ldots,U_n))Min_{1\leq l \leq n}ord_t(\varphi ^{1*}(f_l))>i \frac{p}{q}Min_{1\leq l\leq n}ord_t(\varphi^{*1}(f_l).$$
Ce qui est contradictoire avec la d\'efinition de $\frac{p}{q}$.\endproof

\begin{cor}\textrm{ }\\
Soient $(A,\mathfrak{m})$ un anneau local r\'egulier d'\'egale caract\'eristique z\'ero de dimension $n$ et $I$ un id\'eal $\mathfrak{m}$-primaire alors tout \'el\'ement $x$ de $\overline{I}$ satisfait une relation de d\'ependance int\'egrale sur $I$ de degr\'e $e(I)$ o\`u $e(I)$ d\'esigne la multiplicit\'e de Samuel de $I$.
\end{cor}
\noindent \textit{Preuve :}\\
D'apr\`es les r\'esultats de [N-R], on peut trouver $f_1,\ldots,f_n \in I$ tels que $f_1,\ldots,f_n$ soit une suite r\'eguli\`ere de $A$ et $J=(f_1,\ldots,f_n)$ une r\'eduction minimale de $I$. On a alors $\overline{I}=\overline{J}$. Soit $x \in \overline{I}=\overline{J}$. Si $\hat{A}$ d\'esigne le compl\'et\'e $\mathfrak{m}$-adique de $A$, on a $e(I)=e(I.\hat{A})=e(J)=e(J.\hat{A})$. D'apr\`es le r\'esultat pr\'ec\'edent, $x$ satisfait une relation de d\'ependance int\'egrale de degr\'e $e(I)$ sur $J.\hat{A}$ et donc \`a fortiori sur $I.\hat{A}$. En proc\'edant comme dans [H-I-O] lemme 4.11 p.19 ceci implique que :
$$ (I.\hat{A}+x.\hat{A})^{e(I)}=I.\hat{A}.(I.\hat{A}+x.\hat{A})^{e(I)-1}$$
Comme $\hat{A}$ est fid\'element plat sur $A$, on a :
$$ (I.\hat{A}+x.\hat{A})^{e(I)}\cap A=(I+x.A)^{e(I)} \textrm{ et }I.\hat{A}.(I.\hat{A}+x.\hat{A})^{e(I)-1}\cap A=I.(I+x.A)^{e(I)-1}$$
Par suite : $(I+x.A)^{e(I)}=I.(I+x.A)^{e(I)-1}$ et ceci implique que $x$ satisfait une relation de d\'ependance int\'egrale de degr\'e $e(I)$ sur $I$ toujours en reprenant la preuve de 4.11 de [H-I-O].\endproof 

\section{Preuve de 1.1}
\subsection{Existence du i\`eme exposant de \L ojasiewicz de $I$}
Soient $(A,\mathfrak{m},k)$ un anneau local r\'egulier de dimension $n+1$, d'\'egale caract\'eristique z\'ero, et $I$ un id\'eal $\mathfrak{m}$-primaire. Notons d'abord que le morphisme canonique $A\longrightarrow \hat{A}$ \'etant fid\`element plat, la fonction asymptotique de Samuel est invariante par passage au compl\'et\'e. Ainsi $\overline{v_{I}}(g)=\overline{v_{I.\hat{A}}}(g)$ et $\overline{v_{I.A_H}}(g)=\overline{v_{I.\hat{A_H}}}(g)$ pour tout $g\in A$ et tout $i$ plan $H$. On pourra donc sans restriction supposer que $A$ est complet. Par le th\'eor\`eme de structure de I.S. Cohen c.f. [M], on peut donc supposer que $A=k[[X_0,X_1,\ldots,X_n]]$ avec $Car(k)=0$. Soit $I=(f_1,\ldots,f_m)$  un id\'eal $\mathfrak{m}$ primaire de $A$. D'apr\`es les r\'esultats sur la r\'eduction des id\'eaux de Northcott-Rees [N-R], [H-S], on peut supposer que pour $j>n+1$, $f_j$ est entier sur $J=(f_1,\ldots,f_{n+1})$ i.e. $J$ est une r\'eduction de $I$. Il en d\'ecoule bien \'evidemment que $J.A_H$ est une r\'eduction de $I.A_H$ pour tout $i$-plan $H$. Nous pouvons donc supposer que $I=(f_1,\ldots,f_{n+1})$.
D'autre part, par r\'ecurrence sur la dimension de $A$, il suffit d'\'etablir 1.1 pour les $n$-plans, c'est \`a dire pour une section hyperplane de $A$.\\
Rappelons maintenant que $K=(g_1,\ldots,g_n,h_{n+1})$ est dit une r\'eduction jointe (c.f. [H-S] chap. 17) de $I,\ldots,I,\mathfrak{m}$ ($I$ list\'e $n$ fois), ou encore une r\'eduction jointe de $I^{[n]},\mathfrak{m}$, si et seulement si $g_i \in I$, $h_{n+1}\in \mathfrak{m}$, et $(g_1,\ldots,g_n). I^{n-1}\mathfrak{m}+h_{n+1}.I^n$ est une r\'eduction de $I^n.\mathfrak{m}$. C'est \`a dire s'il existe un entier $k\in \mathbb{N}^*$ tel que :
$$ ((g_1,\ldots,g_n)I^{n-1}\mathfrak{m}+h_{n+1}.I^n).(I^n.\mathfrak{m})^k=I^{n(k+1)}\mathfrak{m}^{k+1}.$$
Ce qui s'\'ecrit encore :
$$ (g_1,\ldots,g_n)I^{n(k+1)-1}\mathfrak{m}^{k+1}+h_{n+1}.I^{n(k+1)}\mathfrak{m}^k=I^{n(k+1)}\mathfrak{m}^{k+1}.$$
Posons :
$$A'=A/(h_{n+1}), \quad J'=(g_1,\ldots,g_n).A' \textrm{ et } I'=I.A/(h_{n+1})=I+(h_{n+1}).A/(h_{n+1}).$$ 
Alors $J'$ est une r\'eduction de $I'$. En effet, si $\mathfrak{m'}$ est l'id\'eal maximal de $A'$, on a :
$$ (*)\hspace{10pt}(g_1,\ldots,g_n)I'^{n(k+1)-1}\mathfrak{m}'^{k+1}=I'^{n(k+1)}.\mathfrak{m}'^{k+1}$$ 
Comme $J'\subset I'$, pour toute valuation discr\`ete $v$ sur $A'$ on a :$J'.A'_v\subset I'.A'_v$ i.e. $v(J')\geq v(I')$. R\'eciproquement par l'\'egalit\'e :
$$(*)\hspace{2pt} v(J')+(n(k+1)-1)v(I')+(k+1)v(\mathfrak{m}')=n(k+1)v(I')+(k+1)v(\mathfrak{m}'),$$
on obtient apr\`es simplification $v(J')=v(I')$. Par cons\'equent $J'\subset I'\subset \overline{J'}$. Il existe maintenant d'apr\`es [H-S], [R-S] des ouverts de Zariski denses $U_1\subset (k^{n+1})^n\simeq (I/\mathfrak{m}.I)^n$ et $V_1\subset k^{n+1}\simeq \mathfrak{m}/\mathfrak{m}^2$ tels que si :
$$g_1=\sum_{1\leq i \leq n+1}\lambda _{i,1}.f_i,\ldots,g_n=\sum_{1\leq i \leq n+1}\lambda _{i,n}.f_i \textrm{ et } h_{n+1}=\sum_{0\leq i \leq n}a_i.X_i$$
satisfont $(\lambda_{i,j}) \in U_1$  et $(a_i)_{0\leq i \leq n}\in V_1$ alors $K=(g_1,\ldots,g_n,h_{n+1})$ est une r\'eduction jointe de $I^{[n]},\mathfrak{m}$. Notons $V_0$ l'ouvert de Zariski de $k^{n+1}$ d\'efini par $X_0\neq 0$. Notant $W_0$ son intersection avec $V_1$, on en d\'eduit qu'il existe un ouvert de Zariski $W_1$ de $k^n$ tel que si   
$$g_1=\sum_{1\leq i \leq n+1}\lambda _{i,1}.f_i,\ldots,g_n=\sum_{1\leq i \leq n+1}\lambda _{i,n}.f_i \textrm{ et } h_{n+1}=X_0-\sum_{1\leq i \leq n}a_i.X_i$$
satisfont $(\lambda _{i,j})\in U_1$ et $(a_i)_{1\leq i \leq n}\in W_1$ alors $(g_1,\ldots,g_n,X_0-\sum_{1\leq i \leq n}a_iX_i)$ est une r\'eduction jointe de $I^{[n]},\mathfrak{m}$. Fixons un $\Lambda =(\lambda _{i,j})$ dans $U_1$ et notons encore $g_1,\ldots,g_n$ les \'el\'ements de $k[[X_0,\ldots,X_n]]$ correspondants. Pour $1\leq i \leq n$, on d\'efinit alors  les \'el\'ements suivants de $k[A][[X_1,\ldots,X_n]]=k[a_1,\ldots,a_n][[X_1,\ldots,X_n]]$ :
$$G_i(A,X_1,\ldots,X_n)=g_i(a_1X_1+\ldots+a_nX_n,X_1,\ldots,X_n)\in k[A][[X_1,\ldots,X_n]].$$
D\'esignons par $\mathfrak{a}\subset k[A][[X_1,\ldots,X_n]]$ l'id\'eal engendr\'e par les $G_i(A,X)$, $1\leq i\leq n$, et soit $\bigtriangleup \subset \mathbb{N}^n$ son diagramme des exposants initiaux au sens de la section pr\'ec\'edente. Pour $a \in k^n$ donn\'e, notons $J'_a\subset k[[X_1,\ldots,X_n]]$ l'id\'eal engendr\'e par les $G_i(a,X)$ apr\`es \'evaluation des coefficients en $a$ et soit $\bigtriangleup _a\subset \mathbb{N}^n$ son diagramme des exposants initiaux. Par construction, pour $a \in W_1$, $J'_a$ est une r\'eduction de $I'_ak[[X_1,\ldots,X_n]]=I.k[[X_0,\ldots,X_n]]/(X_0-\sum_{1\leq i\leq n}a_iX_i)$
On a alors le lemme suivant qui est  une simple transposition des lemmes 7.1 et 7.2 de [B-M].
\begin{lem}\textrm{ }\\
1) $\forall a \in k^n$, $\triangle \leq \triangle _a$.\\
2) Il existe $Q_1,\ldots,Q_l \in k[A]$ tels que $\Gamma=\{a \in k^n/Q_1(a)\times\ldots \times Q_l(a)=0\}$ soit strictement inclus dans $k^n$ et :\\
a) $\forall a \in W_2=k^n-\Gamma $, $\bigtriangleup =\bigtriangleup _a$\\
b) Si $\beta ^1,\ldots,\beta ^l$ d\'esignent les sommets de $\bigtriangleup $, il existe $R_i \in \mathfrak{a}$, $1\leq i \leq l$ tels que :
$$\forall a \in W_2=k^n-\Gamma,\quad \nu (R_i(A,X))=\nu (R_i(a,X))=\beta ^i$$
\end{lem} 
On peut en fait prendre pour $Q_i$ le coefficient du mon\^ome initial de $R_i(X)=\sum _{\alpha \in \mathbb{N}^n}r_{i,\alpha }X^\alpha $, $r_{i,\alpha }\in k[A]$. Notons maintenant $S\subset k[A]$ la partie multiplicativement ferm\'ee engendr\'ee par les $Q_i$, i.e. $S=\{Q_1^{m _1}\times\ldots\times Q_t^{m_t}/(m_1,\ldots,m_t)\in k^t\}$.\\
Soit $\overline{\bigtriangleup }$ (resp. $\overline{\bigtriangleup _a}$) le compl\'ementaire de $\bigtriangleup $ (resp. $\bigtriangleup _a$) dans $\mathbb{N}^n$. $\overline{\bigtriangleup }$ est n\'ecessaire\-ment un ensemble fini puisque dans l'ouvert de Zariski dense $W_1\cap W_2$ de $k^n$ on a : $\overline{\bigtriangleup }=\overline{\bigtriangleup }_a$, et ce dernier ensemble est fini car $\sqrt{J'_a}=(X_1,\ldots,X_n)$. Ainsi pour tout $a \in W=W_1\cap W_2$, les $X^\alpha $, $\alpha \in \overline{\Delta }$, constituent une base du quotient $k[[X_1,\ldots,X_n]]/J'_a$. Il en r\'esulte que pour tout $a \in W$ le morphisme 
$\theta ^{*}_a$:  $k[[Y_1,\ldots,Y_n]]\longrightarrow k[[X_1,\ldots,X_n]]$ d\'efini par $Y_i\rightarrow G_i(a,X)$ fait de $k[[X_1,\ldots,X_n]]$ un $k[[Y_1,\ldots,Y_n]]$ module libre de base les $X^{\alpha}$, $\alpha \in \overline{\bigtriangleup }$. Donc tout \'el\'ement $g$ de $k[[X_1,\ldots,X_n]]$ s'\'ecrit de mani\`ere unique :
$$g=\sum _{\alpha \in\overline{\bigtriangleup }}a_{\alpha}(G_1(a,X),\ldots,G_n(a,X)).X^{\alpha },\quad a_{\alpha }\in k[[Y]].$$ 
Nous avons besoin de d\'ecrire la variation des $a_{\alpha }$ en fonction de $a$ pour \'etablir la constance de $\overline{v_{I.A_H}}(\mathfrak{m}_H)$ sur un ouvert de Zariski. Pour cela notons $R=S^{-1}k[A]$ et soit $\theta ^*$: $R[[Y_1,\ldots,Y_n]]\longrightarrow R[[X_1,\ldots,X_n]]$ qui \`a $Y_i$ fait correspondre $G_i(A,X)$. Nous allons constater que : $$(\bullet)\hspace{5pt} R[[X]] \textrm{ est via } \theta ^*\textrm{ un }R[[Y]] \textrm{ module de type fini engendr\'e par les }X^{\alpha },\hspace{2pt}\alpha \in \overline{\bigtriangleup }.$$
En effet, soit $D \in R[[X]]$. L'algorithme de division formelle de Grauert-Hironaka (c.f section 2) permet d'\'ecrire de mani\`ere unique :
$$D=\sum _{1\leq i \leq t}D_i(X)R_i(X)+\sum _{\alpha \in \overline{\bigtriangleup }}r_{\alpha }X^{\alpha }$$
avec $r_{\alpha} \in R$, $D_i\in R[[X]]$ et $Supp(D_i)+\beta ^i\subset \Delta _i=(\beta ^i+\mathbb{N}^n)-\cup _{k<i}(\beta ^k+\mathbb{N}^n)$. Ce qui se r\'e\'ecrit, en tenant compte du fait que $(R_1,\ldots,R_t)\subset \mathfrak{a}=(G_1,\ldots,G_n)$ en :
$$D=\sum _{1\leq i \leq n}C_i(X)G_i(X)+\sum _{\alpha \in \overline{\bigtriangleup }}r_{\alpha }X^{\alpha }.$$
Le terme $\sum _{\alpha \in \overline{\bigtriangleup }}r_{\alpha }X^{\alpha }$ restant unique puisque son support est inclus dans $\overline{\bigtriangleup }$. Une it\'eration de ce proc\'ed\'e permet de conclure \`a $(\bullet)$. En effet, supposons que pour $k\in \mathbb{N}^*$ nous disposions d'une \'ecriture :
$$(*_k)\hspace{3pt}D=\sum _{\gamma   \in \mathbb{N}^n/ \vert \gamma  \vert=k}C_{\gamma   }G^{\gamma   }+\sum _{\alpha \in \overline{\bigtriangleup }}(\sum _{\gamma   \in \mathbb{N}^n/ \vert\alpha \vert<k}r_{\alpha,\gamma  } G^{\gamma  })X^{\alpha }$$
avec $G^{\gamma   }=G_1(A,X)^{\gamma   _1}\times \ldots \times G_n(A,X)^{\gamma  _n }$ et $\gamma  =(\gamma _1,\ldots,\gamma  _n)$. Divisons par l'algorithme d'Hironaka-Grauert chaque $C_{\gamma }  $ par $R_1,\ldots,R_t$ puis retournant \`a une \'ecriture en les $G_i$, on obtient :
$$C_{\gamma  } =\sum _{1\leq i \leq n}C_{\gamma   ,i}(X)G_i(X)+\sum _{\alpha \in \overline{\bigtriangleup }}a_{\gamma ,\alpha }X^{\alpha }, \quad a_{\gamma ,\alpha }\in R.$$
Reportant dans $(*_k)$, on obtient une \'ecriture :
$$(*_{k+1})\hspace{5pt}D=\sum _{\gamma   \in \mathbb{N}^n/ \vert \gamma  \vert=k+1}C_{\gamma   }G^{\gamma   }+\sum _{\alpha \in \overline{\bigtriangleup }}\left(\sum _{\gamma   \in \mathbb{N}^n/ \vert \gamma   \vert<k+1}r_{\alpha,\gamma  } G^{\gamma }\right)X^{\alpha }.$$
Ceci prouve que d\'esignant par $M\subset R[[X]]$ le $R[[Y]]$ module engendr\'e via $\theta ^*$  par les $X^{\alpha }$, $\alpha \in \overline{\bigtriangleup }$, on a :
$$\forall k\in \mathbb{N}^*,\quad R[[X]]\subset M+\mathfrak{a}^k\subset M+(X)^k.$$
Donc, d'apr\`es le th\'eor\`eme d'intersection de  Krull, on a $R[[X]]=M$ et donc $(\bullet)$. En fait l'unicit\'e dans le th\'eor\`eme de division formel nous dit en plus que $R[[X]]$ est un $R[[Y]]$ module \textit{libre} de type fini via $\theta ^*$ dont une base est constitu\'ee par les $X^{\alpha} $, $\alpha \in \overline{\bigtriangleup }$. Ceci va nous permettre de calculer \textit{globalement} nos polyn\^omes caract\'eristiques et d'obtenir la constance de $\overline{v_{I.A_H}}(\mathfrak{m}_H)$ sur un ouvert de Zariski. En effet, notons $e_1$ le cardinal de $\overline{\bigtriangleup }$ (qui n'est autre que la multiplicit\'e mixte $e(I^{[n]},\mathfrak{m})$ c.f. [H-S]) et soit $i$, $1\leq i \leq n$. On peut \'ecrire pour tout $\alpha \in \overline{\bigtriangleup }$ :
$$X_i.X^{\alpha }=\sum _{\beta \in \overline{\bigtriangleup }}C^i_{\alpha ,\beta }(G_1(A,X),\ldots,G_n(A,X)).X^{\beta }$$
avec $C^i_{\alpha ,\beta }\in R[[Y]]$. Soit alors $M_i$ la matrice carr\'ee d'ordre $e_1$ \`a coefficients dans $R[[Y]]$ :
$$M_i=\left(C^i_{\alpha ,\beta }\right)_{\alpha ,\beta \in \overline{\bigtriangleup }}.$$
Pour chaque $a$ dans l'ouvert de Zariski $W=W_1\cap W_2$, \textit{l'\'evaluation} en $a$, $M_i(a)$ de $M_i$ est la matrice de l'op\'erateur de multiplication par $X_i$ dans le $k[[Y]]$ module libre de type fini $k[[X]]$ via $\theta ^*_a$ d\'etermin\'e par $Y_i\longrightarrow G_i(a,X)$. Notons $P_i(\lambda ,A,Y)=det(\lambda Id_{e_1}-M_i)\in R[[Y]][\lambda ]$ le polyn\^ome caract\'eristique de $M_i$ :
$$P_i(\lambda ,A,Y)=\lambda ^{e_{1}}+\sum _{1 \leq k \leq e_1}r^i_k(A,Y).\lambda ^{e_1-k},\quad r^i_k(A,Y)\in R[[Y]].$$
Ainsi pour chaque $a \in W$, \textit{l'\'evalu\'e} en $a$ de $P_i(\lambda ,A,Y)$ not\'e $P_i(\lambda ,a,Y)$ est le polyn\^ome caract\'eristique d\'esir\'e. Posons pour tout $i,k$, $1\leq i \leq n$, $1\leq k\leq e_1$ : $d^{i}_k=ord_{(Y)}r^{i}_k(A,Y)$ ($ord_{(Y)}( \quad)$ comme \'el\'ement de $R[[Y]]$) et enfin :
$$\overline{v_I}^{(n)}=Min_{1\leq i \leq n}(Min_{1\leq k\leq e_1}\frac{d^i_k}{k}).$$
On peut \'ecrire pour tout $i,k$, $1\leq i \leq n$, $1\leq k \leq e_1$ :
$$r^i_k(A,Y)=\sum _{\gamma \in N^n/\vert \gamma \vert=d^i_k}\frac{r^i_{k,\gamma }(A)}{Q^i_k(A)}Y^{\gamma }+S^i_k(Y)$$
avec $S^i_k(Y)\in (Y)^{d^i_k+1}R[[Y]]$, $r^i_{k,\gamma }(A) \in k[A]$, $Q^i_k(A)\in S$. Consid\'erons alors :
$$V^i_k=\{a \in k^n/\forall \gamma \in \mathbb{N}^n \textrm{ avec } \vert \gamma \vert=d^i_k,\hspace{2pt} r^i_{k,\gamma }(a)=0\}  \textrm{ et } U^i_k=(k^n-V_k^i)\cap W.$$
Par construction chaque ouvert de Zariski $U^i_k$ est non vide. Alors pour toute forme lin\'eaire $H(X_0,X_1,\ldots,X_n)=X_0-\sum _{1\leq i \leq n}a_iX_i$ telle que $a=(a_1,\ldots,a_n)\in U=\cap _{i,k}U^i_k\cap W$, on a :
$$\overline{v_{I.A_H}}(\mathfrak{m}_H)=\overline{v_I}^{(n)}.$$
Ceci au vu du calcul de la fonction asymptotique de Samuel fait \`a la section 3. On obtient ainsi le premier point de 1.1.\endproof
\subsection{Majoration de la multiplicit\'e} 
La preuve relativement simple de l'in\'egalit\'e 2) de 1.1 se fait par r\'ecurrence sur $n$ (dim A=n+1). Elle d\'ecoule d'une g\'en\'eralisation de la loi d'associativit\'e pour les multiplicit\'es que l'on peut trouver dans [N].
\begin{theo}\textrm{([N] Chap.7 Th.18 p. 342)}\\
Soit $(A,\mathfrak{m},k)$ un anneau local noetherien de dimension $s$ et $E$ un $A$-module de type fini. Consid\'erons $a_1,\ldots,a_s$ des \'el\'ements de $A$ engendrant un id\'eal $\mathfrak{m}$-primaire. Alors pour tout $i$, $0\leq i\leq s$, on a :
\begin{eqnarray*}
\lefteqn{e(a_1,\ldots,a_s, E)=}\\
& \sum _{P\in Min(a_1,\ldots,a_i)}e_{A_P}(\phi _P(a_1),\ldots,\phi _P(a_i),E_p)e_{\frac{A}{P}}(\psi _P(a_{i+1}),\ldots,\psi _P(a_s),A/P)
\end{eqnarray*}
o\`u $P$ parcourt l'ensemble des id\'eaux premiers minimaux contenant $(a_1,\ldots,a_i)$ et $\phi _P$, $\psi _P$ d\'esignent respectivement les morphismes canoniques  $\phi _P:$ $A\longrightarrow A_P$ et $\psi _P :$  $A\longrightarrow \frac{A}{P}$.
\end{theo}
Nous prouvons \`a pr\'esent l'in\'egalit\'e 2) de 1.1 par r\'ecurrence sur $n$, $dim A=n+1$. Si $n=0$, il n'y a rien \`a prouver car alors $e(I)=\nu ^{(1)}_I(\mathfrak{m})=ord_{\mathfrak{m}}(I)$. On supposera donc $n>0$ et le r\'esultat \'etabli pour tout anneau local r\'egulier d'\'egale caract\'eristique z\'ero $B$ de dimension $n$. Pr\'esentons $I$ (ou plut\^ot une r\'eduction de $I$) sous la forme $(g_1,\ldots,g_n,g_{n+1})$ o\`u comme pr\'ec\'edemment $K=(g_1,\ldots,g_n,h_{n+1})$ est une r\'eduction jointe de $I^{[n]},\mathfrak{m}$ et $g_{n+1}$ est une combinaison lin\'eaire g\'en\'erique d'un syst\`eme de g\'en\'erateurs de $I$. Ainsi comme nous l'avons vu au paragraphe pr\'ec\'edent on a :
$$g_{n+1}\in \overline{(g_1,\ldots,g_n).\frac{A}{(h_{n+1})}}=\overline{I.\frac{A}{(h_{n+1})}}.$$
Soient $P_1,\ldots,P_l$ les id\'eaux premiers minimaux de $A$ contenant $(g_1,\ldots,g_n)$. Pour tout $j$, $1\leq j\leq l$, on a $dim (A/P_j)=1$. Soit $\widetilde{A/P_j}$ la cl\^oture int\'egrale de $A/P_j$ dans son corps des fractions. $\widetilde{A/P_j}$ est un anneau local noetherien de dimension $1$ et int\'egralement clos, c'est donc un anneau local r\'egulier de dimension 1 et par suite un anneau de valuation discr\`ete. Notons $v_{p_j}$ la valuation d\'efinie par son id\'eal maximal     :
$$v_{P_j}(h)=long(\frac{\widetilde{A/P_j}}{(h)})=Max\{k\in \mathbb{N}/ h \in \mathfrak{m}_{P_j}\}=ord_{\mathfrak{m}_{P_j}}(h).$$
Maintenant pour tout $g \in A$, d\'esignant encore par $g$ l'image de celui-ci via le morphisme naturel $A\longrightarrow A/P_J\longrightarrow \widetilde{A/P_j}$, on a :
$$(*)\hspace{10pt}e_{A/P_j}(g,A/P_j)=long(\frac{A/P_j}{(g)})=long(\frac{\widetilde{A/P_j}}{(g)})=v_{P_j}(g).$$
Maintenant, pour chaque $j$, $1\leq j\leq l$, on a $v_{P_j}(I)=v_{P_j}(g_{n+1})$. D'autre part puisqu'on peut choisir les coefficients de $h_{n+1}=\sum _{i=0}^{n}\lambda _iX_i$ dans un ouvert de Zariski dense de $\mathfrak{m}/\mathfrak{m}^2\simeq k^{n+1}$, quitte \`a restreindre cet ouvert on peut supposer que :
$$\forall j,\quad 1\leq j \leq l,\quad v_{P_j}(\mathfrak{m}_A)=v_{P_j}(h_{n+1}).$$
Ainsi pour tout $j$, $1\leq j\leq l$, $$\frac{v_{P_j}(g_{n+1})}{v_{P_j}(h_{n+1})}=\frac{v_{P_j}(I)}{v_{P_j}(\mathfrak{m}_A)}\leq \nu_{I}^{(n+1)}.$$
En effet, $\nu_{I}^{(n+1)}=Sup \frac{v(I)}{v(\mathfrak{m}_A)}$ car $\nu_{I} ^{(n+1)}=\overline{v_I}(\mathfrak{m}_A)^{-1}=(Min \frac{v(\mathfrak{m}_A)}{v(I)})^{-1}$, les $Sup$ et $Min$ \'etant pris sur l'ensemble des valuations discr\`etes de rang 1 de $A$ (c.f. [H-S]). Nous pouvons \`a pr\'esent appliquer 3.2 avec $i=n$. On a :
$$e(I)=e(g_1,\ldots,g_n,g_{n+1},A)=\sum _{1\leq j\leq l}e(\phi _{P_j}(g_1),\ldots,\phi _{P_j}(g_n),A_{P_j})e(\psi _{P_j}(g_{n+1}),A/P_j).$$
Ce qui s'\'ecrit gr\^ace \`a $(*)$ en :
$$e(I)=e(g_1,\ldots,g_n,g_{n+1},A)=\sum _{1\leq j\leq l}e(\phi _{P_j}(g_1),\ldots,\phi _{P_j}(g_n),A_{P_j})v_{P_j}(g_{n+1})$$
Ecrivant $v_{P_j}(g_{n+1})=v_{P_j}(h_{n+1}).\frac{v_{P_j}(g_{n+1})}{v_{P_j}(h_{n+1})}=v_{P_j}(h_{n+1}).\frac{v_{P_j}(I)}{v_{P_j}(\mathfrak{m}_A)}$ et majorant cette derni\`ere fraction par $\nu _I^{n+1}$, on obtient :
$$e(I)=e(g_1,\ldots,g_n,g_{n+1},A)\leq \left(\sum _{1\leq j\leq l}e(\phi _{P_j}(g_1),\ldots,\phi _{P_j}(g_n),A_{P_j})v_{P_j}(h_{n+1})\right).\nu _I^{(n+1)}.$$
Mais toujours d'apr\`es $3.2$, le terme entre parenth\`eses n'est autre que :
 $$e(g_1,\ldots,g_n,h_{n+1},A)=e(g_1,\ldots,g_n,A/(h_{n+1}).$$
Ainsi : $$e(I)\leq (e(g_1,\ldots,g_n,A/(h_{n+1}))\nu_I ^{(n+1)}.$$
Mais puisque $\overline{(g_1,\ldots,g_n).\frac{A}{(h_{n+1})}}=\overline{I.\frac{A}{(h_{n+1})}}$, on a :
$$e(I,A/(h_{n+1}))=e(g_1,\ldots,g_n,A/(h_{n+1})).$$
Ainsi $e(I)\leq (e(I,A/(h_{n+1})).\nu_I ^{(n+1)}$. Il suffit alors pour conclure d'appliquer l'hypoth\`ese de r\'ecurrence dans l'anneau $B=A/(h_{n+1})$ \`a $I. B$. On a alors $e(I,B)\leq \Pi _{1\leq i \leq n}\nu _I^{(i)}$. Ce qui fournit l'\'egalit\'e 2) de 1.1. \endproof
\begin{rem}\textrm{ }\\
1) Avec les notations ci-dessus le terme $e(I,A/(h_{n+1})=e(g_1,\ldots,g_n,h_{n+1},A)$ n'est autre que la \textrm{multiplicit\'e mixte} $e(I^{[n]},\mathfrak{m},A)$ et nous renvoyons \`a [H-S] chap.17 pour les d\'efinitions et notations.\\
2) On a en fait prouv\'e l'in\'egalit\'e :
$$\frac{e(I^{[n+1-i]},\mathfrak{m}^{[i]},A)}{e(I^{[n-i]},\mathfrak{m}^{[i+1]},A)}\leq \nu _I^{(n+1-i)}.$$
De m\^eme $e(I^{[n+1-i]},\mathfrak{m}^{[i]},A)\leq \Pi _{1\leq j\leq n+1-i} \nu _{I}^{(j)}$.
\end{rem}

\section{Sur les cas d' \'egalit\'e}
Nous cherchons ici comment se caract\'erise les id\'eaux $I$ tels que $e(I)=\Pi _{i=1}^{n+1}\nu _I^{i}$ (Les notations et hypoth\`eses sont celles de 1.1). Pour cela soient $(A,\mathfrak{m},k)$ un anneau local r\'egulier de dimension $n+1$ et $Gr_{\mathfrak{m}}(A)=\bigoplus_{k\in \mathbb{N}}\frac{\mathfrak{m}^k}{\mathfrak{m}^{k+1}}\simeq k[X_0,\ldots,X_n]$. Si $g\in A-(0)$ et si $a=ord_{\mathfrak{m}}(g)$, on appellera forme initiale de $g$ et on notera $In(g)$ la classe de $g$ dans $\frac{\mathfrak{m}^a}{\mathfrak{m}^{a+1}}\subset Gr_{\mathfrak{m}}(A)$. Celle-ci s'identifie \`a un polyn\^ome homog\`ene de degr\'e $a$ de $k[X_0,\ldots,X_n]$, pour tout choix d'un syst\`eme r\'egulier de param\`etres $X_0,\ldots,X_n $ de $A$. Si $I$ est un id\'eal de $A$, on notera $In(I)$ l'id\'eal de $Gr_{\mathfrak{m}}(A)$ engendr\'e par les formes initiales des \'el\'ements $g$ de $I$. Des \'el\'ements $f_1,\ldots,f_m$ de $I$ sont dit une $\mathfrak{m}$-base standard de $I$ si et seulement si $In(g_1),\ldots, In(g_m)$ engendrent $In(I)$. On a alors le r\'esultat suivant.
\begin{prop}\textrm{ }\\
Soient $(A,\mathfrak{m},k)$ un anneau local r\'egulier d'\'egale caract\'eristique z\'ero et $I$ un id\'eal $\mathfrak{m}$-primaire. Consid\'erons les propri\'et\'es suivantes :\\
\begin{itemize}
\item[1)] $e(I)= \nu _I^{(1)}\times\nu ^{(2)}_I\times \ldots.\times \nu _I^{(n)}\times \nu _I^{(n+1)}$
\item[2)] Il existe $b\in \mathbb{N}^*$ et $g_1,\ldots,g_{n+1} \in A$ tels que :\\
a) $In(g_1),\ldots,In(g_{n+1}))$ sont sans z\'eros communs non triviaux dans $\overline{k}^{n+1}$ o\`u $\overline{k}$ est une cl\^oture alg\'ebrique de $k$,\\
b) $\overline{(g_1,\ldots,g_{n+1})}=\overline{I^b}$.
\end{itemize}
Alors $2\Longrightarrow 1$ et si $dim A=2$ alors $1)\Longleftrightarrow 2)$.
\end{prop}

\noindent \textit{Preuve :}\\
Commen\c{c}ons par constater que $2)\Longrightarrow 1)$, ce qui est \'el\'ementaire. Soient $g_1,\ldots,g_{n+1}$ et $b\in \mathbb{N^*}$ satisfaisant $2)$. Posons $a_i=ord_{\mathfrak{m}}(g_i)$ et indexons $g_1,\ldots,g_{n+1}$ de telle sorte que $a_1\leq a_2\leq \ldots a_{n+1}$. Soit $J=(g_1,\ldots,g_{n+1})$. Puisque $\overline{J}=\overline{I^b}$, on a :
$$e(J)=b^{n+1}e(I) \textrm{ et }\nu _J^{(i)}=\nu _{I^b}^{(i)}=b.\nu ^{(i)}_I,\quad 1\leq i\leq n+1$$
Par suite pour obtenir $1)$ pour $I$, il suffit de l'obtenir pour $J$. Pour cela, il nous suffira de constater :
 $$ (*)\hspace{15pt}\forall i,\quad  1\leq i \leq n+1, \quad \nu _J^{(i)}=a_i.$$ En effet, on a d'abord par une propri\'et\'e bien classique de la multiplicit\'e (c.f. [M] Th 14.9 p. 109): $$\Pi _{1\leq i \leq n+1}ord_{\mathfrak{m}}(g_i)=\Pi _{1\leq i \leq n+1}a_i\leq e(J).$$
L'\'egalit\'e s'obtient par $(*)$ en utilisant la majoration $2)$ de 1.1.\\
Prouvons $(*)$.  Pour cela, il suffit de prouver que si $g_1,\ldots,g_{n+1}$ sont $n+1$ \'el\'ements de $A$ satisfaisant $2a)$ et index\'es selon $ord_{\mathfrak{m}}()$ croissant alors : $\nu ^{(n+1)}_J=a_{n+1}$ o\`u $J=(g_1,\ldots,g_{n+1})$. En effet supposons cette affirmation  prouv\'ee en toute g\'en\'eralit\'e. Alors soit $h(X_0,\ldots,X_n)=X_0-\sum_{I=1}^na_iX_i$ une forme lin\'eaire suffisamment g\'en\'erale pour que $In(g_1),\ldots,In(g_n),h$ soient sans z\'eros communs non triviaux dans $\overline{k}^{n+1}$. D\'esignons par $g'_j$ la classe de $g_j$ dans $A'=A/(h)$, $1\leq j \leq n+1$. Notons $J_1'=(g'_1,\ldots,g'_n)$ et $J'=(g'_1,\ldots,g'_n,g'_{n+1})$. Alors $In(g'_1),\ldots, In(g'_n)$ sont sans z\'eros communs non triviaux dans $\overline{k}^n$ (ce ne sont autres que $In(g_1)(\sum_{1\leq i \leq n}a_iX_i,X_1\ldots,X_n)$ $,\ldots,In(g_n)(\sum_{1\leq i \leq n}a_iX_i,X_1\ldots,X_n))$. Par cons\'equent si notre assertion est prouv\'ee $\nu _{J'_1}^{(n)}=a_n$. Maintenant soit $\mathfrak{m}'$ le maximal de $A'$. On a $g'_{n+1}\in \mathfrak{m}'^{a_{n+1}}\subset \mathfrak{m}'^{a_n}$ et $\mathfrak{m}'^{a_n}\subset \overline{J'_1}$ car $\nu _{J'_1}^{(n)}=a_n$. Par suite $g'_{n+1}\in \overline{J'_1}$ et donc $\overline{J'}=\overline{J'_1}$ et $\nu _{J'}^{(n)}= \nu _{J'_1}^{(n)}=a_n$ i.e $\nu _{J}^{(n)}=a_n$. Il suffit alors de r\'ep\'eter l'op\'eration. \\
 Il ne nous reste plus qu'a prouver  que sous $2a)$ on a $\nu _J^{(n+1)}=a_{n+1}$. L'hypoth\`ese 2)a) fait que $g_1,\ldots,g_{n+1}$ est une suite r\'eguli\`ere et une $\mathfrak{m}$-base standard de $J=(g_1,\ldots,g_{n+1})$ (c.f. [H-I-O] 13.10 p 96). Soit $l\in \mathbb{N}^*$ assez grand pour que $\mathfrak{m}^l$ soit inclus dans $J$. A fortiori : $\mathfrak{m}^{la_{n+1}}\subset J$. Puisque $g_1,\ldots g_{n+1}$ est une base standard de $J$, on peut ecrire pour tout $r \in \mathfrak{m}^{la_{n+1}}$  (c.f. [H.I.O] 13.7 p. 91 ):
$$r=\sum_{i=1}^{n+1}q_ig_i, \textrm{ avec } ord_{\mathfrak{m}}(q_i) \geq la_{n+1}-a_i \geq (l-1)a_{n+1}.$$
Par suite : 
$$(\mathfrak{m}^{a_{n+1}}+J)^l\subset (\mathfrak{m}^{a_{n+1}}+J)^{l-1}.J\subset (\mathfrak{m}^{a_{n+1}}+J)^l$$
Ainsi $J$ est une r\'eduction de $J+\mathfrak{m}^{a_{n+1}}$ et ces deux id\'eaux ont donc m\^eme cl\^oture int\'egrale i.e. $\overline{J}=\overline{J+\mathfrak{m}^{a_{n+1}}}$. En particulier : $\mathfrak{m}^{a_{n+1}}\subset \overline{J}$. De ce dernier fait, on d\'eduit que $\nu _J^{(n+1)}\leq \nu _{\mathfrak{m}^{a_{n+1}}}^{(n+1)}=a_{n+1}$. L'in\'egalit\'e oppos\'ee 
s'obtient en consid\'erant un arc (non trivial) $\varphi ^{*}:$ $A\longrightarrow \overline{k}[[t]]$ tel que $\varphi ^*(g_1)=\ldots=\varphi ^*(g_n)=0$ (ce qui est possible car $Dim A/(g_1,\ldots,g_n)=1$). Comme $In(g_1),\ldots,In(g_{n+1})$ sont sans z\'eros communs non triviaux dans $\overline{k}^{n+1}$, on a $ord (\varphi ^{*}(g_{n+1}))=a_{n+1}ord(\varphi ^{*}(\mathfrak{m}))$. Ceci fournit une valuation $v$ telle que $v(I)/v(\mathfrak{m})\geq a_{n+1}$ et donc $\nu _j^{(n+1)}\geq a_{n+1}$.\\
Montrons maintenant que si $dim A=2$, alors $1\Longrightarrow 2)$. Posons $\nu ^{(2)}_I=\frac{a_2}{b}$ et $\nu ^{(1)}_I=ord_{\mathfrak{m}}(I)=a_1$. Soit $g_1 \in I$ tel que $ord_{\mathfrak{m}}(g_1)=a_1$. On peut trouver un syst\`eme r\'egulier de param\`etres de $A$, $X_0,X_1$,  tel que $In(g_1)(X_0,0)\neq 0$ (ceci quitte \`a effectuer un changement \guillemotleft lin\'eaire \guillemotright\quad  de syst\`eme r\'egulier de param\`etres). Donc $ord_{\mathfrak{m'}}(g_1.A/(X_1))=a_1$ o\`u $\mathfrak{m'}$ est le maximal de $A/(X_1)$. Puisque $\nu ^{(2)}_I=a_2/b$, on a : $X_1^{a_2}\in \overline{I^b}$. Par cons\'equent : $J=(X_1^{a_2},g_1^b)\subset \overline{I^b}$. Mais $e(J)= a_2.b.e(X_1,g_1)$ et $e(X_1,g_1)=a_1$ car $ord_{\mathfrak{m'}}(g_1.A/(X_1))=a_1$. Ainsi : $e(J)=a_2.b.a_1=b^2(\frac{a_2}{b}\times a_1)=b^2e(I)=e(I^b)$. Par cons\'equent, par un c\'el\`ebre r\'esultat de  D. Rees (c.f. [H-S] Th 11.3.1 p.222) on a $\overline{J}=\overline{I^b}$, et donc 2).\endproof
\begin{rem}\textrm{ }\\
Nous ignorons si en g\'en\'eral si on a \'equivalence entre les conditions $1)$ et $2)$. En fait, soient $A$ comme ci-dessus avec $dim A=n+1$, $n\geq 2$, et $I$ un id\'eal satisfaisant $1)$. Comme nous avons vu en 4.3 que :
$$e(I)\leq e(I^{[n]},\mathfrak{m})\nu _I^{(n+1)}\leq \Pi _{1\leq k\leq n+1}\nu _I^{(k)}.$$
On a :
$$e(I^{[n]},\mathfrak{m})=\Pi _{1\leq k\leq n}\nu _I^{(k)} \textit{ et }e(I^{[n+1-i]},\mathfrak{m}^{[i]})=\Pi _{1\leq k\leq n+1-i}\nu _I^{(k)} .$$
Ceci conduit a une caract\'erisation de $1)$ du type $2)$ mais seulement apr\`es section hyperplane g\'en\'erique. Par exemple si $dim A=3$, notant $\nu _I^{(k)}=\frac{a_i}{b}$. Si $(X_0,X_1,X_2)$ est un syst\`eme r\'egulier de param\`etres de $A$ suffisamment g\'en\'eral pour que $(g_1,g_2,X_2)$ soit une une r\'eduction jointe de $(I^{[2]},\mathfrak{m})$, on obtient comme pr\'ec\'edemment que $\overline{I^b}=\overline{(g_1^b,g_2^b,X_2^{a_3})}$. Posons : $K=I^b.A/(X_3)$. Alors $K$ satisfait fait $1)$, et donc par le cas $dim A=2$, il existe $g'_1,g'_2 \in \overline{K}$ tels que $\overline{K}=\overline{(g'_1,g'_2)}$ et $In(g'_1)(X_0,X_1), In(g'_2)(X_0,X_1)$ sont sans z\'eros communs non triviaux dans $\overline{k}^2$. Il se pose alors la question de rel\'evement suivante :\\
Existe t-il $h_1,h_2\in A$ tels que :\\
- $h_1,h_2 \in \overline{I^b}$ et  classe de $h_i$ dans $A/(X_2)$ \'egale $g'_i$,\\
- $In(h_i)(X_0,X_1,0)=In(g'_i)(X_0,X_1)$.\\
Nous ignorons en g\'en\'eral la r\'eponse \`a de telles questions.
\end{rem}
\vspace{10pt}

\centerline{ \textit{BIBLIOGRAPHIE }}

\vspace{10pt}

\noindent [A-H-V] J.M. Aroca-H. Hironaka-J.L. Vicente, \textit{The theory of the maximal contact}, Mem. Mat. Inst. Jorge Juan n°29, Consejo Superior de Investigaciones Cientificas, Madrid, 1975.

\noindent [B1] N. Bourbaki, \textit{Alg\`ebre Chapitre 1 \`a 3}, nouvelle \'edition Hermann 1970.

\noindent [B2] N. Bourbaki, \textit{Alg\`ebre commutative chap.1 \`a 7}, Hermann 1961. 

\noindent [B3] N. Bourbaki, \textit{Alg\`ebre commutative chap. 10}, Masson 1998.

\noindent [B-M] E. Bierstone, P.D. Milman, \textit{Relations among analytic functions I}, Ann. Inst. Fourier $\mathbf{37}$ (1),(1987) 187-239. 

\noindent[E] D. Eisenbud, \textit{Commutative algebra with a view toward algebraic geometry}, Graduate Texts in Mathematics $n^{°}$ 150, Springer-Verlag.

\noindent[G] H. Grauert, \textit{Uber die deformation isolierter Singularitaten analytischer Mengen}, Invent. Math. $\mathbf{15}$ (1972), 171-198.

\noindent [H-I-O] M. Herman, S. Ikeda, U. Orbanz, \textit{Equimultiplicity and blowing up, an algebraic study with an appendix by B. MOONEN}, Springer 1988.

\noindent [H-S] C. Huneke-I. Swanson, \textit{Integral Closure of ideals, Rings, and Modules}, London Mathematical Society Lecture Note Series $n^{°}$ 336, Cambridge University Press 2006.

\noindent [La] J.P. Lafon, \textit{Alg\`ebre commutative, langages g\'eom\'etrique et alg\'ebrique}, Collection enseignement des Sciences n° 24, Hermann.

\noindent [Lo] S. \L ojasiewicz, \textit{Ensembles semi-analytiques}, Pub. Math. I.H.E.S. 1964.

\noindent [L-T] M. Lejeune-B. Teissier, \textit{Cl\^oture int\'egrale des id\'eaux et \'equisingularit\'e}, S\'emi\-naire Lejeune-Teissier, Centre de Math\'ematiques de l'\'ecole polytechnique 1974, Publications Universit\'e Scientifique et M\'edicale de Grenoble.

\noindent [M] H. Matsumura, \textit{Commutative Ring theory}, Cambridge studies in mathematics 8, 1986.

\noindent [N] D.G. Northcott, \textit{Lessons on Rings, Modules and Multiplicities}, Cambridge University Press 1968.

\noindent [N.R] D.G. Northcott-D. Rees, \textit{Reduction of ideals in local rings}, Proc. Cambridge Phil. Soc., $\mathbf{50}$ 2 (1954), 145-158.

\noindent [R1] D. Rees, \textit{Lectures on the asymptotic theory of ideals}, London Mathematical lecture notes series 113, Cambridge University Press.

\noindent [R2] D. Rees, \textit{Multiplicities, Hilbert functions an degree functions.} In Commutative algebra: Durham 1981,  London Math. Soc. Lecture Note Ser., $\mathbf{72}$, Cambridge-New-York, Cambridge University Press, 1982, 170-178.

\noindent [R-S] D. Rees-J. Sally, \textit{General elements and joint reductions}, Michigan Math. J., $\mathbf{35}$ (1988), 241-254.

\noindent [T1] B. Teissier, \textit{Cycles \'evanescents, sections planes et conditions de Whitney}, Ast\'erisque 7 et 8, (1973), 285-362.

\noindent[T2] B. Teissier, \textit{Vari\'et\'es Polaires I, Invariants polaires des singularit\'es d'hypersurfaces}, Invent. Math. $\mathbf{40}$ (3), (1977), 267-292.

\noindent [T3] B. Teissier, \textit{Sept compl\'ements au s\'eminaire Lejeune-Teissier}, \`a para\^{\i}tre dans Ann. Fac. Sci. Toulouse.

\noindent [To] J.C. Tougeron, \textit{Id\'eaux de fonctions diff\'erentiables}, Ergebnisse der Mathematik und ihrer Grenzgebiete Band 71, Springer-Verlag 1972.

 Michel HICKEL,\\ 
Universit\'e Bordeaux 1, I.M.B.\\
Equipe d'Analyse et G\'eom\'etrie\\
et I.U.T. Bordeaux 1 d\'epartement Informatique\\ 
33405 Talence Cedex, France\\
email : hickel@math.u-bordeaux1.fr 
\end{document}